\documentclass[a4paper,oneside,english]{amsart}
\usepackage[T1]{fontenc}
\usepackage[latin1]{inputenc}
\usepackage{setspace}
\usepackage{amssymb}

\makeatletter

\providecommand{\tabularnewline}{\\}

 \theoremstyle{plain}
\newtheorem{thm}{Theorem}[section]
  \theoremstyle{definition}
  \newtheorem{defn}[thm]{Definition}
  \theoremstyle{remark}
  \newtheorem{rem}[thm]{Remark}

\usepackage{babel}
\makeatother
\begin{document}

\title{An elementary sieve}

\date{August 28th, 2007}

\author{Damián Gulich$^{A,B}$\\
Gustavo Funes$^{A,B}$\\
Leopoldo Garavaglia$^{C}$\\
Beatriz Ruiz$^{D}$\\
Mario Garavaglia$^{A,B}$}

\maketitle
\begin{singlespace}
\begin{center}
{\small $^{A}$Departamento de Física, Facultad de Ciencias Exactas,
Universidad Nacional de La Plata, Argentina}
\par\end{center}{\small \par}

\begin{center}
{\small $^{B}$Laboratorio de Procesamiento Láser, Centro de investigaciones
Ópticas, La Plata, Argentina}
\par\end{center}{\small \par}

\begin{center}
{\small $^{C}$Aranjuez, Spain}\\
{\small $^{D}$Centro de investigaciones Ópticas, La Plata, Argentina}
\par\end{center}{\small \par}
\end{singlespace}

\begin{singlespace}

\email{\emph{\small E-mail addresses:} {\small Damián Gulich: dgulich@ciop.unlp.edu.ar;
Gustavo Funes: gfunes@ciop.unlp.edu.ar; Mario Garavaglia: garavagliam@ciop.unlp.edu.ar}}
\end{singlespace}

\begin{abstract}
In this paper we review the properties of families of numbers of the
form $6n\pm1$, with $n$ integer (in which there are all prime numbers
greater than 3 and other compound numbers with particular properties)
to later use them in a new sieve that allows the separation of numbers
$n$ that generate primes from those that only generate compounds.
In principle, this can be used to find the amount of prime numbers
up to a given number $h$; this means, $\pi\left(h\right)$.
\end{abstract}

\keywords{\emph{Key words and phrases:} Prime numbers, sieve, prime counting
function.}

\section{Introduction}

An old problem in mathematics is the way to compute the amount of
prime numbers less or equal to a given value $h$ \cite[p. 347]{Dickson}.
This function is known as $\pi\left(h\right)$ \cite{Weisstein}.
The preeminent method for such task since the 3rd century BC was the
sieve of Eratosthenes, and there were no great advances on the subject
until the work of Gauss in 1863 \cite[p. 352]{Dickson}; which also
allowed more advances \cite{Hardy,Deleglise} along with the growth
of calculation power in the 20th century. Recent implementations \cite{Deleglise}
require diverse relatively elaborated computational strategies.

In this paper we will study this problem defining a new sieve whose
properties permit an elementary study of $\pi\left(h\right)$, and
also the possibility of finding its value on a given interval.

\section{About the form $6n+1$}

We begin with a very well known result:

\begin{thm}
\label{thm:Modulos}Every prime number of absolute value greater than
3 can be written in the form $6n+1$ or $6n-1$.
\end{thm}
\begin{proof}
Let's see the equivalences modulo 6. Suppose $q$ prime.

1) If $q\underset{6}{=}0$ $\Rightarrow$$q=6n$$\Rightarrow$$6\mid q$,
ABS.

2) If $q\underset{6}{=}1$ $\Rightarrow$$q=6n+1$, which is not impossible
since $7=6+1$ is a prime.

3) If $q\underset{6}{=}2$ $\Rightarrow$$q=6n+2=2\left(3n+1\right)$$\Rightarrow$$2\mid q$,
ABS.

4) If $q\underset{6}{=}3$ $\Rightarrow$$q=6n+3=3\left(2n+1\right)$$\Rightarrow$$3\mid q$,
ABS.

5) If $q\underset{6}{=}4$ $\Rightarrow$$q=6n+4=2\left(3n+2\right)$$\Rightarrow$$2\mid q$,
ABS.

6) If $q\underset{6}{=}5$ $\Rightarrow$$q=6n+5=6n+6-1=6\left(n+1\right)-1$,
which is not impossible since with $n=1$ this gives 11, a prime.
\end{proof}
\begin{defn}
The $\alpha$ class of integer numbers \cite{Garavaglia} is the set
\begin{equation}
\alpha=\left\{ x\in\mathbb{Z}/x=6n+1,n\in\mathbb{Z}\right\} \label{eq:6n+1}\end{equation}

\end{defn}
~

\begin{defn}
The $\beta$ class of integer numbers \cite{Garavaglia} is the set

\begin{equation}
\beta=\left\{ x\in\mathbb{Z}/x=6n-1,n\in\mathbb{Z}\right\} \label{eq:6n-1}\end{equation}

\end{defn}
The different values of relations in \eqref{eq:6n+1} and \eqref{eq:6n-1}
are shown in table \ref{tab:Valores-para-6n-1}.

\begin{center}
\begin{table}
\begin{centering}
\begin{tabular}{|c|c|c|}
\hline 
$n$&
$\beta_{n}=6n-1$&
$\alpha_{n}=6n+1$\tabularnewline
\hline
\hline 
-5&
-31 {*}&
-29 {*}\tabularnewline
\hline 
-4&
-25&
-23{*}\tabularnewline
\hline 
-3&
-19 {*}&
-17 {*}\tabularnewline
\hline 
-2&
-13 {*}&
-11 {*}\tabularnewline
\hline 
-1&
-7 {*}&
-5 {*}\tabularnewline
\hline 
0&
-1&
1\tabularnewline
\hline 
1&
5 {*}&
7 {*}\tabularnewline
\hline 
2&
11 {*}&
13 {*}\tabularnewline
\hline 
3&
17 {*}&
19 {*}\tabularnewline
\hline 
4&
23 {*}&
25\tabularnewline
\hline 
5&
29 {*}&
31 {*}\tabularnewline
\hline
\end{tabular}
\par\end{centering}

\caption{\label{tab:Valores-para-6n-1}Values of $6n-1$ and $6n+1$ for several
$n$. With ({*}) we mark prime numbers.}
\end{table}

\par\end{center}

Strictly speaking, a {}``complete'' list of all prime numbers of
absolute value greater than 3 ($\left\{ \ldots,-7,-5,5,7,\ldots\right\} $)
is the list of primes from both classes.

We now state a property given in \cite{Garavaglia} where product
rules are proved as a theorem: 

\begin{center}
\begin{table}[!h]
\begin{tabular}{|c||c|c|}
\hline 
$\times$&
$\alpha$&
$\beta$\tabularnewline
\hline
\hline 
$\alpha$&
$\alpha$&
$\beta$\tabularnewline
\hline 
$\beta$&
$\beta$&
$\alpha$\tabularnewline
\hline
\end{tabular}

\caption{\label{tab:Tabla-de-productos.}Table of products.}
\end{table}

\par\end{center}

\begin{thm}
Every prime number of absolute value greater than 3 (except for the
sign) is generated by $6n+1$, with $n$ integer.
\end{thm}
\begin{proof}
We must prove the equivalence (except for the sign) between both families
given in Theorem \ref{thm:Modulos}. 

Let be $f_{\alpha}\left(n\right)=6n+1$ and $f_{\beta}\left(n\right)=6n-1$,
we mus now prove that $f_{\alpha}\left(-n\right)=-f_{\beta}\left(n\right)$.
Indeed:\[
f_{\alpha}\left(-n\right)=6\left(-n\right)+1=-6n+1=-\left(6n-1\right)=-f_{\beta}\left(n\right)\]

\end{proof}
\begin{defn}
We define the set of integer numbers $G_{\alpha}$:\begin{equation}
G_{\alpha}=\left\{ g\in\mathbb{Z}/6g+1\mbox{ is a prime}\right\} \label{eq:G+}\end{equation}

This means, $G_{\alpha}$ is the set of \emph{all} numbers that (except
for the sign) generate \emph{all} primes of absolute value greater
than 3 by the relationship \eqref{eq:6n+1}.
\end{defn}

\section{The sieve}

\begin{defn}
Let $A$ be an infinite matrix whose element $a\left(i,j\right)$
\footnote{Coordinates are in the Cartesian sense.%
} is \begin{equation}
a\left(i,j\right)=i+j\left(6i+1\right)\label{eq:a_ij}\end{equation}
where $i,j\,\in\mathbb{Z}$.
\end{defn}
Note that numbers on the axis also match this representation.

\[
\begin{array}{ccccccccccc}
 &  &  &  &  & \vdots\\
 & -96 & -71 & -46 & -21 & 4 & 29 & 54 & 79 & 104\\
 & -73 & -54 & -35 & -16 & 3 & 22 & 41 & 60\\
 & -50 & -37 & -24 & -11 & 2 & 15 & 28\\
 & -27 & -20 & -13 & -6 & 1 & 8\\
\cdots & -4 & -3 & -2 & -1 & 0 & 1 & 2 & 3 & 4 & \cdots\\
 & 19 & 14 & 9 & 4 & -1\\
 & 42 & 31 & 20 &  & -2\\
 & 65 & 48 &  &  & -3\\
 & 88 &  &  &  & -4\\
 &  &  &  &  & \vdots\end{array}\]

\subsection{Properties}

\begin{thm}
\label{thm:A-es-sim=E9trica.}$A$ is symmetrical.
\end{thm}
\begin{proof}
A simple expansion shows that \[
a\left(i,j\right)=i+j\left(6i+1\right)=i+6ij+j=j+i\left(6j+1\right)=a\left(j,i\right)\]

\end{proof}
\begin{defn}
Let $\widetilde{A}$ be the set of unrepeated elements of $A$ excluding
the axes (the elements of the form $a\left(i,0\right)$ y $a\left(0,j\right)$).
\end{defn}

\subsubsection{About the signs of the elements of $\tilde{A}$}

Four quadrants can be seen:
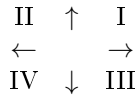
\[
\begin{array}{ccc}
\mbox{II} & \uparrow & \mbox{I}\\
\leftarrow &  & \rightarrow\\
\mbox{IV} & \downarrow & \mbox{III}\end{array}\]

What happens to the signs of the elements of $\tilde{A}$ from each
quadrant?

Because of Theorem \ref{thm:A-es-sim=E9trica.}, we should only focus
on the sign of elements of $\tilde{A}$ originally from quadrants
I, II, and IV.

\begin{itemize}
\item In quadrant I ($i\geq1$, $j\geq1$) all elements are positive
\item In quadrant II ($i\leq-1$, $j\geq1$)\\
It's easy to see that $j\left(6i+1\right)\leq0$, then $\tilde{a}\left(i,j\right)\leq0$
$\forall i,j$. 
\item In quadrant IV ($i\leq-1$, $j\leq1$)\\
$i+j\left(6i+1\right)=i+j+6ij$. $\left(i+j\right)\leq-1$ y $ij\geq\left|i+j\right|\geq1$.
Then, the sign is positive.
\end{itemize}
\begin{thm}
\label{thm:Los-elementos-de-A-tilde-NO}The elements of $\widetilde{A}$
DO NOT generate prime numbers.
\end{thm}
\begin{proof}
$\tilde{a}\left(i,j\right)=i+j\left(6i+1\right)$ with $i\neq0$ y
$j\neq0$. If we put this into \eqref{eq:6n+1} and suppose $p$ prime,
then\[
p=6\tilde{a}\left(i,j\right)+1=6\left(i+j\left(6i+1\right)\right)+1=6i+6j\left(6i+1\right)+1=\left(6i+1\right)+6j\left(6i+1\right)=\left(6i+1\right)\left(6j+1\right)\]
 but since $i$ and $j$ are different from zero, then $p$ would
be a compound, ABS.
\end{proof}
\begin{rem}
\label{rem:Posibilidades}According to the signs of $i$ and $j$,
$6\tilde{a}\left(i,j\right)+1$ sweeps (except for the sign) all possibilities:
\end{rem}
\begin{enumerate}
\item If $i>0$ y $j>0$, the generated number is of the form $\alpha\cdot\alpha$.
\item If $i<0$ y $j<0$, the generated number is of the form $\beta\cdot\beta$.
\item If $i$ y $j$ have opposite signs, the generated number is of the
form $\alpha\cdot\beta$.
\end{enumerate}
See Table \ref{tab:Tabla-de-productos.} for properties of products
of $\alpha$'s and $\beta$'s.

\section{Leopoldo's Theorem}

\begin{thm}
\label{thm:Teorema-de-Leopoldo}(Leopoldo's Theorem) $G_{\alpha}=\mathbb{Z}-\widetilde{A}$
\end{thm}
\begin{proof}
Every prime number of absolute value greater than 3 is either $\alpha$
or $\beta$, and products between those classes of equivalence are
closed on themselves. Because of Theorem \ref{thm:Los-elementos-de-A-tilde-NO}
and Remark \ref{rem:Posibilidades}, we know that $\widetilde{A}$
generates all possible $\alpha$ and $\beta$ compound numbers, and
so the elements $n$ of $\mathbb{Z}\notin\widetilde{A}$ generate
prime numbers by $6n+1$.
\end{proof}

\section{Conclusions}

Numbers in $\widetilde{A}$ are easy to generate. However, this values
are repeated. If an ordering of its elements is possible, then comparing
them with $\mathbb{Z}$ will give a list of the values in $G_{\alpha}$.
The existence of the $\widetilde{A}$-like elements has been known
for a long time \cite[p. 356]{Dickson}, but have been considered
impractical to find primes. In this sieve representation, the number
of primes between two given values $h_{1}=6c_{1}-1$ and $h_{2}=6c_{2}+1$
where $c_{1}$ and $c_{2}$ and both are greater than zero, would
be \begin{equation}
\pi\left(h_{2}\right)-\pi\left(h_{1}\right)=2\left(c_{2}-c_{1}\right)-\Lambda\left(c_{1},c_{2}\right)+\xi\label{eq:delta}\end{equation}
 where $\Lambda\left(c_{1},c_{2}\right)$ would be a procedure that
would count all non repeated values of $\widetilde{A}$ and $\xi$
is a fitting factor. A way of doing this will be the subject of a
future paper. Given the value of $\pi\left(h_{1}\right)$, $\pi\left(h_{2}\right)$
would be computable with \eqref{eq:delta}.

\section{Acknowledgments}

Damián Gulich and Gustavo Funes are financially supported by a student
fellowship from the INNOVATEC Foundation, Argentina.

Damián Gulich and Gustavo Funes thank Dr. Mario Garavaglia for involving
them in this line of research.

Damián Gulich dedicates this paper to MMB and to his new niece Antonia
Monti.

\end{document}